\documentclass{article} 

\usepackage{amssymb,amsmath}
\usepackage{graphicx}

\begin{document}\centerline{\bf How to Verify the Riemann Hypothesis}\vskip .3in

\centerline{ M.L. Glasser}\vskip .2in
\centerline{Department of Physics}
\centerline{Clarkson University}
\centerline{Potsdam, Ny 13699-5820 (USA)}\vskip .4in
\centerline{\bf ABSTRACT}\vskip .1in

\begin{quote}
It is proposed that  the validity, or not,  of the Riemann Hypothesis might be established on the basis of the integral
$$\int\frac{\xi(2s)}{\xi(s)}ds$$
where
$$\xi(s)=(s-1)\pi^{-s/2}\Gamma(1+s/2)\zeta(s).$$
\end{quote}

\newpage

 In his famous Monograph of 1859 {\it  On the Number of Primes less than a certain Magnitude}[1] Riemann introduced the entire function
 $$\xi(s)=(s-1)\pi^{-s/2 }\Gamma(1+s/2)\zeta(s)$$
 The related function
  $f(\sigma+i\tau)=\xi[2(\sigma+i\tau)]/\xi(\sigma+i\tau)$ is meromorphic with all its poles to the left of $\sigma=1$ and on $\sigma=1/2$ if Riemann's hypothesis, that all the non-trivial zeroes of $\zeta(\sigma+i\tau)$  lie on the line $\sigma=1/2$, is true.  Let
  $$I(\epsilon)=\int_{1/2+\epsilon-i\infty}^{1/2+\epsilon+i\infty}\frac{\xi(2t)}{\xi(t)}\frac{dt}{2\pi i} =\frac{1}{\pi}\int_0^{\infty}Re[f[(1/2+\epsilon+it)]dt,  \quad \epsilon>0.$$
  Then,  since $f$ is known to vanish at infinity, by Cauchy's theorem
$$\int_{1/2+\epsilon-i\infty}^{1/2+\epsilon+i\infty} f(t)\frac{dt}{2\pi i}$$
$$=\frac{1}{\pi}\int_0^{\infty}Re[f(\frac{3}{2}+it)]dt+  Re[S(\epsilon)].$$
where $S(\epsilon)$ Is the sum of the residues of any poles $ \sigma+i\tau$ with  $1/2<\sigma<1$, $\tau>0$.  If the Riemann hypothesis is false and  $\zeta(2s)$, $\zeta(s)$ have no common zeros,  $S(\epsilon)\ne 0$.. Thus $I(\epsilon)$ being independent of $\epsilon$ would suggest that Riemann's conjecture   is true.

Fig.1  is a plot of $I(\epsilon)$ vs $\epsilon$

\begin{figure}[htbp]
\centering
\includegraphics[width=6cm]{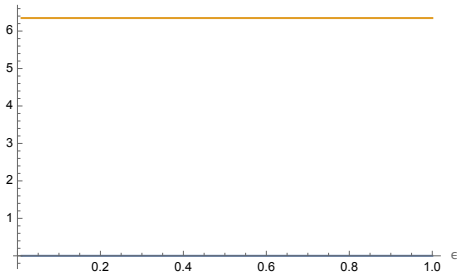}
\caption{$I(\epsilon)$ vs $\epsilon$.}
    \label{fig1}
    \end{figure}
    
\noindent
Unfortunately, Fig.1 is based on only  six decimal places, but
the problem would be settled definitively could the exact value of $I(\epsilon)$ be known explicitly, which seems  a daunting task, although it has been possible to show that  [2]
$$\int_{-\infty}^{\infty} \xi(\frac{1}{2}+\epsilon+it)dt=\frac{\pi^{1/4}}{\sqrt{32}}\Gamma(1/4)\left[\frac{\Gamma^8(1/4)}{32\pi^4}-3\right] $$
is independent of $\epsilon$.
\vskip .1in
\newpage
\noindent
{\bf Acknowledgement}\vskip .1in
The author thanks Dr. M. Milgram for comments.\vskip .2in

\centerline{\bf REFERENCES}\vskip .1in

\noindent
[1] H.M. Edwards, {\it Riemann's Zeta Function}  New York: Academic Press, ISBN 0-12-232750-0, Zbl 0315.10035

\noindent
[2] M.L. Glasser,  SCIENTIA
Series A: Mathematical Sciences, Vol. 31 (2021), 15-23(ISSN 0716-8446)

\noindent
[3] E.C. Titchmarsh {\it Theory of the Riemann Zeta Function},[ Oxford UP (1953)] (1.2.11)

\end{document}